# Linearization of multivariate discrete difference operators
*Yunting Iris Gao*

# Contents

# 1 Abstract

In 2023 in (3), Uwe finds the explicit form of the map which is which is settled in $\mathbb{Z}^N$ of finite functional degree and discusses how to compute its usual degree w.r.t to the derivative in the linear form, i.e. the product of ones formed by its orthogonal basis, and also introduces the notion of functional degree. It is the linear combination of the product of binomials. And this makes a big progress in the development of integer-valued functions. And this inspires the author to discuss this form in the separate form, i.e. the triple of finite set. In (18), Hrycaj discusses these operators in more abstract form. In this paper, we unifies two approaches. Seriously, we prove that the functional degree of integer-valued maps on integers which is computed with respect to its difference operators associated with its standard basis is exactly the same with the ones computed with respect to arbitrary mix discrete difference operators and we call it as multivariate difference operators. Since the functional degree agrees with the usual degree of integer-valued maps, we actually achieved the linearization of multivariate difference operators

# 2 Main Results

The mains results are stated in this section and they explained the process of linearization and how we catch the notion of the functional degree. Essentially

$$\text{count}(\mathcal{P}) = \text{fdeg}(f(\mathcal{P})) = \text{fdeg}^{\mathbb{Z}^N}(f) = \deg(f(\mathcal{P})) = \text{fdeg}^{\mathbb{A}}(f)$$

The symbol deg is denoted as the usual degree. The second and third equations are discussed in the series of papers (3), (2) and (4) and the last equation is somewhat discussed in (18) and (19) for for arbitrarily polynomial maps.

**Theorem 2.1.** . The counting number of a $\mathbb{Z}$−polyfract $\mathcal{P}$ is determined by the number of taking discrete difference operators of its associated function. Let a triple $(\mathcal{P}, \mathcal{A}, \mathcal{N})$ be a $\mathbb{Z}$−polyfract, i.e some finite $\mathcal{A} \subset \mathbb{Z}$ and some finite $\mathcal{N} \subset \mathbb{N}^N$. Then

$$\text{count}(\mathcal{P}) = \text{fdeg}(f(\mathcal{P}))$$

**Theorem 2.2.** Let $f : \mathbb{Z}^N \to \mathbb{Z}$ be a map of finite functional degree. Then

$$\text{fdeg}^{\mathbb{A}}(f) = \text{fdeg}^{\mathbb{Z}^N}(f)$$

**Theorem 2.3.** Let $h : \mathbb{Z}^N \to \mathbb{Z}$ be a multivariate polynomial such that its corresponding $\mathbb{Z}$−polyfract $\mathcal{P}(\mathcal{A}, \mathcal{N})$ satisfying $\mathcal{A} = \{1\} \subset \mathbb{Z}$ and $\mathcal{N} \subset \mathbb{N}^N_{\geq 1}$ with $\text{card}(\mathcal{N}) = 1$. i.e. $h = f(\mathcal{P})$. Denoting $\mathcal{N} = \{\underline{m}\} \subset \mathbb{N}^N_{\geq 3}$. Let $\underline{n} \in \mathbb{N}^N_{\geq 2}$ such $n_i < m_i$ for each $1 \leq i \leq N$. Then for each $\underline{x} = (x_1, \cdots, x_N) \in \mathbb{Z}^N$:

$$\sum_{0 \leq p_1 \leq n_1; 0 \leq p_2 \leq n_2; \cdots ; 0 \leq p_N \leq n_N} (-1)^{p_1 + \cdots + p_N} h(\sum_{i=1}^N (n_i - p_i + x_i) e_i) = \prod_{l=1}^N \binom{x_l}{m_l - n_l}$$

For most recent progress in integer-valued polynomials and polynomials functions, we recommend (20).

# 3　Preliminaries

Assume $\mathbb{A}, \mathbb{B}$ are commutative groups We recall basic results of arithmetic rules of shift operators and discrete difference operators $[s], \Delta_a$ in $\text{End}(\mathbb{B}^{\mathbb{A}})$ is listed.

**Theorem 3.1.** Assume $\mathbb{A}$ be finitely generated with a generating set $S \subset \mathbb{A}$, i.e $\mathbb{A} = \langle S \rangle_{\mathbb{Z}}$

(a) $\Delta_{-s} = -[-s] \cdot \Delta_s$

(b) $\Delta_{s_1+s_2} = [s_1] \cdot \Delta_{s_2} + \Delta_{s_1}$

(c) $[s]^k = [ks]$ denoting $[s]^k = \underbrace{[s] \cdot [s] \cdot [s] \cdots [s]}_{k\text{-times}}$

(d) (Commutative Law)
$$\Delta_a \cdot \Delta_b = \Delta_a \cdot \Delta_b$$
$$[s_1] \cdot [s_2] = [s_2] \cdot [s_1]$$
$$\Delta_a \cdot [s] = [s] \cdot \Delta_a$$

(e) (Distributive Law)
$$[s] \cdot (\Delta_a + \Delta_b) = [s] \cdot \Delta_a + [s] \cdot \Delta_b \tag{0.1}$$
$$\Delta_a \cdot ([s_1] + [s_2]) = \Delta_a \cdot [s_1] + \Delta_a \cdot [s_2] \tag{0.2}$$

(f) (Additive property)
$$[a + b] = [a] \cdot [b] \tag{0.3}$$

*Proof.* These can be verified by direct computation. ∎

**Theorem 3.2.** Assume $\mathbb{A}$ is commutative. Let $a, b \in \mathbb{A}$. Then $[a]\Delta_b = \Delta_b[a], \Delta_b\Delta_a = \Delta_b\Delta_a$ and $[a][b] = [b][a]$.

*Proof.* These can be verified by direct computation. ∎

We need the following lemma for our reconstruction in the third section.

**Lemma 3.3.** Denote
$$a = r_1 e_1 + r_2 e_2 + \cdots + r_N e_N$$

Then
$$\Delta_a = \sum_{i=1}^{N} \left\{ \left( \prod_{l=1}^{i-1} [e_l]^{r_l} \right) \cdot \left( \sum_{p=0}^{r_i-1} [e_i]^p \right) \right\} \Delta_{e_i}$$

*Proof.* This is a simple calculation from Theorem 3.1 and Theorem 3.2. ∎

We also need the following well-known combinatorial results and however its complete proof is often omitted in most standard textbooks. For references in combinatorics see (21).

**Theorem 3.4** (Multipolynomial Theorem).
$$(x_{1,1} + x_{1,2} + x_{1,3} + \cdots + x_{1,m_1})(x_{2,1} + x_{2,2} + x_{2,3} + \cdots + x_{2,m_2}) \cdots (x_{n,1} + x_{n,2} + x_{n,3} + \cdots + x_{n,m_n})$$
$$= \sum_{p_1=1}^{m_1} \sum_{p_2=1}^{m_2} \cdots \sum_{p_n=1}^{m_n} x_{1,p_1} x_{2,p_2} \cdots x_{n,p_n}$$
$$= \sum_{(p_1, p_2, \cdots, p_n) \in \prod_{i=1}^{n}[m_i]} \left\{ \prod_{i=1}^{n} x_{i,p_i} \right\}$$

Since this Theorem is common, we present its proof in the appendix

# 4 Reconstruction of multivariate discrete difference operators from their generators

In this section, we study multivariate discrete difference operators $\Delta_{\underline{a}} := \Delta_{a_1}\Delta_{a_2}\cdots\Delta_{a_d}$ in $\text{End}(\mathbb{Z}^{\mathbb{Z}^N})$ in the goal-orientated way and we use properties listed in the first section to reconstruct multivariate difference operators to two variations. Essentially the first one is the expansion of its basis. And the second is the expansion in various combination of elements from its basis. All of results in $\text{End}(\mathbb{Z}^{\mathbb{Z}^N})$ can be deduced to $\text{End}(\mathbb{B}^{\mathbb{Z}^N})$

## 4.1 The first variation

$\mathbb{Z}^N$ is a $N-$dimensional finitely generated $\mathbb{Z}-$module and Denote $\underline{e} = (e_1, e_2, \cdots, e_N)$. We are going to find the explicit form $F(\underline{a}, \underline{q}, \underline{e})$ and the condition $C(\underline{a}, \underline{q}, \underline{e})$ for the equation

$$\begin{cases} \Delta_{\underline{a}} = \sum_{C(\underline{a},\underline{q},\underline{e})} F^1(\underline{a}, \underline{q}, \underline{e}) \Delta_{\underline{e}}^{\underline{q}} \\ \qquad \underline{a} \in (\mathbb{Z}^N)^d \\ \qquad \underline{q} \in (\mathbb{Z}^+)^N \end{cases} \quad \text{(Equation 0.3a)}$$

The main technique used in the process of finding explicit form is combinations and relabeling the indices in the multi polynomial products. Denote $\underline{a} = (a_1, a_2, \cdots, a_d)$ and $\underline{q} = (q_1, q_2, \cdots, q_N)$. Choose a $d \times N$ matrix in $\mathbb{N}^{d \times N}$

$$R = (r_{t,m})_{1 \leq t \leq d, 1 \leq m \leq N}$$

Write

$$\begin{aligned} a_1 &= r_{1,1}e_1 + r_{1,2}e_2 + \cdots + r_{1,N}e_N \\ a_2 &= r_{2,1}e_1 + r_{2,2}e_2 + \cdots + r_{2,N}e_N \\ &\vdots \\ a_d &= r_{d,1}e_1 + r_{d,2}e_2 + \cdots + r_{d,N}e_N \end{aligned}$$

For each $1 \leq t \leq d$ and $1 \leq m \leq N$, by Lemma 3.3, denote

$$\alpha_{t,m} = \left\{\left(\prod_{l=1}^{m-1}[e_l]^{r_{t,l}}\right) \cdot \left(\sum_{p=0}^{r_{t,m}-1}[e_m]^p\right)\right\} \in \text{End}(\mathbb{Z}^{\mathbb{Z}^N}) \tag{0.4}$$

Then:

$$\begin{aligned} \Delta_{a_1} &= \alpha_{1,1}\Delta_{e_1} + \alpha_{1,2}\Delta_{e_2} + \cdots + \alpha_{1,N}\Delta_{e_N} \\ \Delta_{a_2} &= \alpha_{2,1}\Delta_{e_1} + \alpha_{2,2}\Delta_{e_2} + \cdots + \alpha_{2,N}\Delta_{e_N} \\ &\vdots \\ \Delta_{a_d} &= \alpha_{d,1}\Delta_{e_1} + \alpha_{d,2}\Delta_{e_2} + \cdots + \alpha_{d,N}\Delta_{e_N} \end{aligned} \tag{0.5}$$

By Theorem 3.4: to form each term we need to choose exactly one term from row $z$ and for the column term we use $\beta_z$ to relate with the row index $z$. Then we have

$$\prod_{z=1}^{d} \alpha_{(z,\beta_z)} \tag{0.6}$$

After observing in each row the index of $\alpha_{t,m}e_m$: we have another product of terms in the summation in $F^1$

$$\Delta_{e_{\beta_1}}\Delta_{e_{\beta_2}}\cdots\Delta_{e_{\beta_d}}$$

Now replacing $(t, m) = (z, \beta_z)$ in the identity 0.4 and substituting it to the product 0.6, we have

$$\prod_{z=1}^{d}\left\{\left(\prod_{l=1}^{\beta_z-1}[e_l]^{r_{(z,\beta_z)}}\right) \cdot \left(\sum_{p=0}^{r_{(z,\beta_z)}-1}[e_z]^p\right)\right\}$$

Now since there are $N$ columns in the R.H.S of 0.5, one possible condition of $C(\underline{a}, \underline{q}, \underline{e})$ is

$$\mathfrak{X}_1 = \{1 \leq \beta_1 \leq N; 1 \leq \beta_2 \leq N; \cdots; 1 \leq \beta_d \leq N\} \tag{Restriction 1}$$

Now we set the following conditions among $\underline{q} = (q_1, \cdots, q_N)$ for each $1 \leq m \leq N$

$$\mathfrak{X}_2 = \{q_m \text{ is multiplicities of } m \text{ among the sequence } (\beta_1, \cdots, \beta_d)\} \tag{Restriction 2}$$

so that

$$\Delta_{e_{\beta_1}} \Delta_{e_{\beta_2}} \cdots \Delta_{e_{\beta_d}} = \Delta_{e_1}^{q_1} \cdots \Delta_{e_N}^{q_N} = \Delta_{\underline{e}}^{\underline{q}}$$

and this can be easily verified by intuition. Now in summary, we have the following result after denoting $[N] = \{1, 2, 3, \cdots, N\}$

$$\Delta_{\underline{a}} = \sum_{C(\underline{a},\underline{q},\underline{e})} \prod_{z=1}^{d} \left\{ \left( \prod_{l=1}^{\beta_z - 1} [e_l]^{r_{(z,\beta_z)}} \right) \cdot \left( \sum_{p=0}^{r_{(z,\beta_z)} - 1} [e_z]^p \right) \right\} \Delta_{\underline{e}}^{\underline{q}} \tag{0.7}$$

where the condition $C(\underline{a}, \underline{q}, \underline{e}) = \mathfrak{X}_1 \bigcup \cup_{m=1}^{N} \mathfrak{X}_2$ We actually show the following statement

**Theorem 4.1.** Denote standard basis of $\mathbb{Z}$–module $\mathbb{Z}^N$ by $\{e_1, e_2, \cdots, e_N\}$. Denote $\underline{e} = (e_1, e_2, \cdots, e_N)$. Denote $\underline{a} = (a_1, \cdots, a_d) \in (\mathbb{Z}^N)^d$. Then

$$\Delta_{\underline{a}} = \sum_{C(\underline{a},\underline{q},\underline{e})} \left\{ \prod_{z=1}^{d} \sum_{p=0}^{r_{(z,\beta_z)} - 1} [\omega_{z,p}] \right\} \Delta_{\underline{e}}^{\underline{q}} \tag{0.8}$$

where for each $1 \leq z \leq d$ and $0 \leq p \leq r_{(z,\beta_z)} - 1$

$$\omega_{z,p} = pe_z + \sum_{l=1}^{\beta_z - 1} r_{(z,\beta_z)} e_l$$

*Proof.* A direct calculation yields the result:

$$\prod_{z=1}^{d} \left\{ \left( \prod_{l=1}^{\beta_z - 1} [e_l]^{r_{(z,\beta_z)}} \right) \cdot \left( \sum_{p=0}^{r_{(z,\beta_z)} - 1} [e_z]^p \right) \right\} = \prod_{z=1}^{d} \left\{ [\sum_{l=1}^{\beta_z - 1} r_{(z,\beta_z)} e_l] \cdot \left( \sum_{p=0}^{r_{(z,\beta_z)} - 1} [e_z]^p \right) \right\} \tag{0.9}$$

$$= \prod_{z=1}^{d} \left\{ [\sum_{l=1}^{\beta_z - 1} r_{(z,\beta_z)} e_l] \cdot \left( \sum_{p=0}^{r_{(z,\beta_z)} - 1} [pe_z] \right) \right\} \tag{0.10}$$

$$= \prod_{z=1}^{d} \left\{ \left( \sum_{p=0}^{r_{(z,\beta_z)} - 1} [pe_z] [\sum_{l=1}^{\beta_z - 1} r_{(z,\beta_z)} e_l] \right) \right\} \tag{0.11}$$

$$= \prod_{z=1}^{d} \left\{ \sum_{p=0}^{r_{(z,\beta_z)} - 1} [pe_z + \sum_{l=1}^{\beta_z - 1} r_{(z,\beta_z)} e_l] \right\} \tag{0.12}$$

where 0.9 and 0.12 are both by the additive property 0.3. ∎

### 4.2 The second variation

The idea to find this variation is almost the same with the first one. But we still present the process. The equations for the second variation is

$$\begin{cases} \Delta_{\underline{a}} = \sum_{C(\underline{a},\underline{e})} F^2(\underline{a},\underline{e}) \Delta_{\underline{e}} \\ \underline{a} \in (\mathbb{Z}^N)^d \end{cases} \tag{Equation 0.12b}$$

Some steps are omitted for finding this variation. To use the idea of partitions, we construct the following matrix

$$R = (r_{t,k})_{1 \leq t \leq d, 1 \leq k \leq N} \in \mathbb{N}^{d \times N}$$

By repeatedly applying Lemma 3.3 to $a_1, a_2, \cdots, a_d$,

$$\Delta_{a_1} = \left\{\left(\prod_{l=1}^{0}[e_l]^{r_{1,l}}\right)\cdot\left(\sum_{p=0}^{r_{1,1}-1}[e_1]^p\right)\right\}\Delta_{e_1} + \left\{\left(\prod_{l=1}^{1}[e_l]^{r_{1,l}}\right)\cdot\left(\sum_{p=0}^{r_{1,2}-1}[e_2]^p\right)\right\}\Delta_{e_2} + \cdots + \left\{\left(\prod_{l=1}^{N-1}[e_l]^{r_{1,l}}\right)\cdot\left(\sum_{p=0}^{r_{1,N}-1}[e_N]^p\right)\right\}\Delta_{e_N}$$

$$\Delta_{a_2} = \left\{\left(\prod_{l=1}^{0}[e_l]^{r_{2,l}}\right)\cdot\left(\sum_{p=0}^{r_{2,1}-1}[e_1]^p\right)\right\}\Delta_{e_1} + \left\{\left(\prod_{l=1}^{1}[e_l]^{r_{2,l}}\right)\cdot\left(\sum_{p=0}^{r_{2,2}-1}[e_2]^p\right)\right\}\Delta_{e_2} + \cdots + \left\{\left(\prod_{l=1}^{N-1}[e_l]^{r_{2,l}}\right)\cdot\left(\sum_{p=0}^{r_{2,N}-1}[e_N]^p\right)\right\}\Delta_{e_N}$$

$$\vdots$$

$$\Delta_{a_d} = \left\{\left(\prod_{l=1}^{0}[e_l]^{r_{d,l}}\right)\cdot\left(\sum_{p=0}^{r_{d,1}-1}[e_1]^p\right)\right\}\Delta_{e_1} + \left\{\left(\prod_{l=1}^{1}[e_l]^{r_{d,l}}\right)\cdot\left(\sum_{p=0}^{r_{d,2}-1}[e_2]^p\right)\right\}\Delta_{e_2} + \cdots + \left\{\left(\prod_{l=1}^{N-1}[e_l]^{r_{d,l}}\right)\cdot\left(\sum_{p=0}^{r_{d,N}-1}[e_N]^p\right)\right\}\Delta_{e_N}$$

Denoting for each $1 \leq t \leq d$ and each $1 \leq m \leq N$

$$\alpha_{t,m} = \left\{\left(\prod_{t=1}^{m-1}[e_l]^{r_{t,l}}\right)\cdot\left(\sum_{p=0}^{r_{t,m}-1}[e_m]^p\right)\right\}$$

Then

$$\Delta_{\alpha_1} = \alpha_{1,1}\Delta_{e_1} + \alpha_{1,2}\Delta_{e_2} + \cdots + \alpha_{1,N}\Delta_{e_N}$$
$$\Delta_{a_2} = \alpha_{2,1}\Delta_{e_1} + \alpha_{2,2}\Delta_{e_2} + \cdots + \alpha_{2,N}\Delta_{e_N}$$
$$\vdots$$
$$\Delta_{a_d} = \alpha_{d,1}\Delta_{e_1} + \alpha_{d,2}\Delta_{e_2} + \cdots + \alpha_{d,N}\Delta_{e_N}$$

We rename the second index of $\alpha_{t,m}$ for $1 \leq t \leq d$ by $k_1, k_2, \cdots, k_d$, by Theorem 3.4:

$$\Delta_{a_1}\Delta_{a_2}\cdots\Delta_{a_d} = \sum_{k_1=1}^{N}\sum_{k_2=1}^{N}\cdots\sum_{k_d=1}^{N}\alpha_{1,k_1}\Delta_{e_{k_2}}\alpha_{2,k_2}\Delta_{e_{k_2}}\cdots\alpha_{d,k_d}\Delta_{e_{k_d}}$$

Now collecting $\alpha_{t,m}$ and discrete difference operators separately in monomials, denoting $[N] := \{1,2,3,\cdots,N\}$:

$$\Delta_{a_1}\Delta_{a_2}\cdots\Delta_{a_d} = \sum_{(k_1,k_2,\cdots,k_d)\in[N]^d}\alpha_{1,k_1}\alpha_{2,k_2}\cdots\alpha_{d,k_d}\Delta_{e_{k_1}}\Delta_{e_{k_2}}\cdots\Delta_{e_{k_d}}$$

$$= \sum_{(k_1,k_2,\cdots,k_d)\in[N]^d}\left(\prod_{t=1}^{d}\alpha_{t,k_t}\right)\Delta_{e_{k_1}}\Delta_{e_{k_2}}\cdots\Delta_{e_{k_d}}$$

Substituting $\alpha_{t,m}$ into it,

$$\Delta_{a_1}\Delta_{a_2}\cdots\Delta_{a_d} = \sum_{(k_1,k_2,\cdots,k_d)\in[N]^d}\prod_{t=1}^{d}\left\{\left(\prod_{l=1}^{k_t-1}[e_l]^{r_{t,l}}\right)\cdot\left(\sum_{p=0}^{r_{t,k_t}-1}[e_{k_t}]^p\right)\right\}\Delta_{e_{k_1}}\Delta_{e_{k_2}}\cdots\Delta_{e_{k_d}}$$

We actually have the following statement

**Theorem 4.2.** Denote the standard basis of $\mathbb{Z}$−module $\mathbb{Z}^N$ by $\{e_1, e_2, \cdots, e_N\}$. Denote $\underline{e} = (e_1, e_2, \cdots, e_N)$. Denote $\underline{a} \in (\mathbb{Z}^N)^d$. Then

$$\Delta_{\underline{a}} = \sum_{(k_1,\cdots,k_d)\in[N]^d}\prod_{t=1}^{d}\left\{\sum_{p=0}^{r_{t,k_t}-1}[\omega_{p,t}]\right\}\Delta_{e_{k_1}}\Delta_{e_{k_2}}\cdots\Delta_{e_d}$$

where $\omega_{p,t} = pe_{k_t} + \sum_{l=1}^{k_t-1}r_{t,l}e_l$ for each $1 \leq t \leq d$ and $0 \leq p \leq r_{t,k_t} - 1$

*Proof.* This is proved by the direct computation and repeatedly using the additive property

$$\prod_{t=1}^{d}\left\{\left(\prod_{l=1}^{k_t-1}[e_l]^{r_{t,l}}\right)\cdot\left(\sum_{p=0}^{r_{t,k_t}-1}[e_{k_t}]^p\right)\right\} = \prod_{t=1}^{d}\left\{\left(\prod_{l=1}^{k_t-1}[r_{t,l}e_l]\right)\cdot\left(\sum_{p=0}^{r_{t,k_t}-1}[e_{k_t}]^p\right)\right\}$$

$$= \prod_{t=1}^{d}\left\{\left[\sum_{l=1}^{k_t-1}r_{t,l}e_l\right]\cdot\left(\sum_{p=0}^{r_{t,k_t}-1}[pe_{k_t}]\right)\right\}$$

$$= \prod_{l=1}^{d}\left\{\sum_{p=0}^{r_{t,k_t}-1}[pe_{k_t}][\sum_{l=1}^{k_t-1}r_{t,l}e_l]\right\}$$

$$= \prod_{t=1}^{d}\left\{\sum_{p=0}^{r_{t,k_t}-1}[pe_{k_t}+\sum_{l=1}^{k_t-1}r_{t,l}e_l]\right\}$$

∎

## 5    Formulas of multivariate discrete difference operators of cyclic groups

The formula found in this section is actually used for the proof of our main result. But since it plays the most important role among theorems used in it. It is discussed in the separate section. Basically we want to find the product $\Delta_{a_1}\cdots\Delta_{a_d}$ of discrete difference operators in $\text{End}(\mathbb{B}^{\mathbb{A}})$ when $\mathbb{A}$ is a cyclic group, i.e there is single element in its generating set. Since a cyclic group is just a simple case of finitely generated group, the result can be deduced from each of the variation in the previous section by carefully thinking about the index under summation. Intuitively we proved that the multivariate discrete difference operators is the product of a discrete difference operator and the $d$-fold of $\Delta_s$ when assuming $\mathbb{A}=\langle s\rangle_{\mathbb{Z}}$.

**Theorem 5.1.** Let $d\in\mathbb{Z}^+$. Let $(a_1,a_2,\cdots,a_d)\in\mathbb{A}^d$. Denote $\mathbb{A}=\langle s\rangle_{\mathbb{Z}}$. Then

$$\Delta_{\underline{a}}=T\Delta_s^d$$

for some $T\in\text{End}(\mathbb{B}^{\mathbb{A}})$

*Proof.* We finished the proof by direct computing using Theorem 3.1. To use it, we need to write by $\mathbb{A}=\langle s\rangle_{\mathbb{Z}}$ for some $s\in\mathbb{A}$

$$a_1=r_1s;\ a_2=r_2s;\ \cdots;a_d=r_ds$$

Because of (a) in Theorem 3.1, without loss of generality, we can assume $r_1,\cdots,r_d\in\mathbb{N}_{\geq 1}$. First of all, we are trying to find $\Delta_{ks}$: applying (b) in Theorem 3.1 to the case $\langle(k-1)s,s\rangle_{\mathbb{Z}}\subset\mathbb{A}$ and inductively

$$\Delta_{rs}=\Delta_{(r-1)s+s}$$
$$=[(r-1)s]\cdot\Delta_s+\Delta_{(r-1)s}$$
$$=[(r-1)s]\cdot\Delta_s+\Delta_{(r-2)s+s}$$
$$=[(r-1)s]\cdot\Delta_s+[(r-2)s]\cdot\Delta_s+\Delta_{(r-2)s}$$
$$=[(r-1)s]\cdot\Delta_s+[(r-2)s]\cdot\Delta_s+\cdots+[s]\cdot\Delta_s$$
$$=\left(\sum_{p=1}^{r-1}[ps]\right)\cdot\Delta_s$$

Applying this formula to the different cases $r_1,r_2,\cdots,r_d$ and denoting $T_i=\sum_{p=1}^{r_i-1}[ps]$ for $1\leq i\leq d$ Then

$$\Delta_{\underline{a}}=\Delta_{a_1}\cdots\Delta_{a_d}$$
$$=\prod_{l=1}^{d}\Delta_{a_l}$$
$$=\prod_{l=1}^{d}\{T_l\Delta_s\}$$
$$=\left\{\prod_{l=1}^{d}T_l\right\}\left\{\prod_{i=1}^{d}\Delta_s\right\}$$
$$=\left\{\prod_{l=1}^{d}T_l\right\}\Delta_s^d$$

where the last step is by Theorem 3.2 and the checking step for $T_l \in \text{End}(\mathbb{B}^{\mathbb{A}})$ is omitted for each $1 \leq l \leq d$. After denoting $T := \prod_{l=1}^{d} T_l$, we finished the proof. ∎

Actually by the deduction from variations or applying Theorem 3.4 with Lemma 3.3, we can find

$$T = \sum_{(p_1,\cdots,p_d)\in\prod_{i=1}^{d}[r_i-1]} \left\{\prod_{i=1}^{d} [s]^{p_i}\right\} \tag{0.13}$$

The technical details of the production process is for each $1 \leq i \leq d$:

$$\Delta_{r_i s} = \left(\sum_{p_i=1}^{r_i-1} [p_i s]\right) \cdot \Delta_s^d$$

Then for $\underline{a} = (a_1, \cdots, a_d)$

$$\Delta_{\underline{a}} = \prod_{i=1}^{d} \Delta a_i$$
$$= \prod_{i=1}^{d} \Delta_{r_i s}$$
$$= \prod_{i=1}^{d} \left(\sum_{p_i=1}^{r_i-1} [p_i s]\right) \Delta_s^d$$
$$= \left\{\prod_{i=1}^{d} \left(\sum_{p_i=1}^{r_i-1} [p_i s]\right)\right\} \Delta_s^d$$

where the last step is by the fact $\text{End}(\mathbb{Z}^{\mathbb{Z}^N})$ is a commutative ring. Applying Theorem 3.4: for each $1 \leq i \leq d$: $[r_i-1] = \{1, 2, \cdots, r_i - 1\}$

$$\prod_{i=1}^{d}\left(\sum_{p_i=1}^{r_i-1}[p_i s]\right) = \sum_{(p_1,\cdots,p_d)\in\prod_{i=1}^{d}[r_i-1]} \prod_{i=1}^{d}[p_i s] = \sum_{(p_1,\cdots,p_d)\in\prod_{i=1}^{d}[r_i-1]} \left\{\prod_{i=1}^{d}[s]^{p_i}\right\}$$

# 6 Linearization of multivariate difference operators

To write a statement formally, keeping the spirit of $\mathbb{Z}$−polyfract in (3), we introduce the notion called the $\mathbb{B}$−polyfract. Basically it can be understood as the sub-collection of $\{\binom{x_1}{n_1}\cdots\binom{x_N}{n_N}b_{\underline{n}} \mid b_{\underline{n}} \in \mathbb{B}, \underline{n} \in \mathbb{N}^N\}$ and hence a $\mathbb{B}$−polyfract can be denoted as a triple

$$(\mathcal{P}, \mathcal{A}, \mathcal{N})$$

where $\mathcal{A} \subset \mathbb{B}$ with $\text{card}(\mathcal{A}) < \infty$ and $\mathcal{N} \subset \mathbb{N}^N$ with $\text{card}(\mathcal{N}) < \infty$ and sometimes denoted by $\mathcal{P}(\mathcal{A}, \mathcal{N})$. We define the counting number of $\mathcal{P}$ as $\max\{|\underline{n}| \mid \underline{n} \in \mathcal{N}\}$. Denoted by $\text{count}(\mathcal{P})$, i.e. $\text{count}(\mathcal{P}) := \max\{|\underline{n}| \mid \underline{n} \in \mathcal{N}\}$. In this notation, $\mathbb{Z}$−polyfract is just a special case of $\mathbb{B}$−polyfract when $\mathbb{B} = \mathbb{Z}$. We denote $f(\mathcal{P})$ as its associated polynomial function and we still use the definition of the functional degree from (3). When there is only one monomial, i.e. the sub-collection is singleton, the $\mathbb{B}$−polyfract triple is $(\mathcal{P}, \{b\}, \{\underline{n}\})$. Roughly speaking, the functional degree is the least number of taking discrete difference operators for vanishing this function minus one or equivalent the maximal number of taking discrete difference operators such that this function is nozero constant or it is not zero which is evaluated at zero. Here we denote $\text{bi} : \mathbb{Z} \times \mathbb{Z} \to \mathbb{Z}$ by $\text{bi}(x,y) = \binom{x}{y}$. Then there are some basical properties about the value of this map which are collected from (3)

**Lemma 6.1.**
- $\text{bi}(x,y) \equiv 0$ if $0 \leq x < y$ and $y \geq 1$
- $\text{bi}(x,y) \equiv 1$ if $y = 0$
- $\text{bi}(x,y) = 0$ if $x = y$
- $\text{bi}(x,y) \equiv 0$ if $y < 0$

**Lemma 6.2.** Let $\mathbb{A} = \langle s \rangle_{\mathbb{Z}}$ be a finitely generated group with one single element. Let $\mathbb{B}$ be a commutative group. Let $d \in \mathbb{Z}^+$. Then the following conditions are equivalent

(a) For each $(h_1, h_2, \cdots, h_d) \in \mathbb{A}^d$: $\Delta_{h_1}\Delta_{h_2}\cdots\Delta_{h_d} \equiv 0$ in $\text{End}(\mathbb{B}^{\mathbb{A}})$.

(b) $\Delta_s^d \equiv 0$ in $\text{End}(\mathbb{B}^{\mathbb{A}})$

*Proof.* The implication $(a) \Rightarrow (b)$ is obvious and just set $h_1 = h_2 = \cdots = h_d = h$. We now prove its reverse statement and arbitrarily choose $d$ elements (possibly pair-wise distinct) $h_1, h_2, \cdots, h_d \in \mathbb{A}$. By $\mathbb{A} = \langle s \rangle_{\mathbb{Z}}$, we choose $d$ integers $r_1, r_2, \cdots, r_d \in \mathbb{Z}^+$ such that $h_1 = r_1 s$, $h_2 = r_2 s$, $\cdots$, $h_d = r_d s$. Denoting

$$T := \sum_{(p_1, \cdots, p_d) \in \prod_{i=1}^d [r_i - 1]} \left\{ \prod_{i=1}^d [s]^{p_i} \right\}$$

for some $T \in \text{End}(\mathbb{B}^{\mathbb{A}})$. Then by 0.13, we have

$$\Delta_{h_1}\Delta_{h_2}\cdots\Delta_{h_d} = T\Delta_s^d \tag{0.14}$$

Since $s \in \langle s \rangle_{\mathbb{Z}} = \mathbb{A}$, applying (b) to R.H.S of the identity 0.14, we have $\Delta_{h_1}\Delta_{h_2}\cdots\Delta_{h_d} = 0$ in $\text{End}(\mathbb{B}^{\mathbb{A}})$ ∎

**Lemma 6.3** ((19), pg.2). Let $f: \mathbb{A} \to \mathbb{B}$ be a map of commutative groups satisfying the following polynomial map with degree less than $k$

- ♣ $f = \sum_{l=0}^{k-1} f_l$, where $f_i$ is a monomial of degree $i$ for each $i$ and $f_0$ constant.

- ★ For every $g_1, g_2, \cdots, g_s \in \mathbb{A}$, there are elements $c_{\underline{n}} \in \mathbb{B}$ and $\underline{n} \in \mathbb{N}^N$ with $|\underline{n}| < k$ such that for every $\eta_1, \cdots, \eta_s \in \mathbb{Z}$

$$f(\eta_1 + \cdots + \eta_s) = \sum_{|\underline{n}| < k} \eta_1^{n_1} \cdots \eta_N^{n_N} c_{\underline{n}}$$

- ♦ For every $a, b \in \mathbb{A}$ there are elements $c_{\underline{n}_1}, \cdots, c_{\underline{n}_k} \in \mathbb{B}$ such that for every $\eta \in \mathbb{Z}$

$$f(a + \eta b) = \sum_{l=1}^k \eta^l c_{\underline{n}_l}$$

Then the following conditions are equivalent

(a) $\Delta_{h_1}\cdots\Delta_{h_k} f = 0$ for each $h_1, \cdots, h_k \in \mathbb{A}^k$

(b) $\Delta_h^k f = 0$ for each $h \in \mathbb{A}$

**Theorem 6.4.** [(3), pg.10] Choose $\underline{n} = (n_1, \cdots, n_N) \in \mathbb{N}^N$. Consider a polyfract polynomial function with single monomial $f: \mathbb{Z}^N \to \mathbb{Z}$, i.e. $f(\underline{x}) = \binom{x_1}{n_1}\cdots\binom{x_N}{n_N}$. Choose $\underline{m} \in \mathbb{Z}^N$ such that $m_i \leq n_i$ for each $i$. Then

$$\Delta^{\underline{m}} f(\underline{x}) = \Delta^{\underline{m}} \prod_{l=1}^N \binom{x_l}{n_l} = \prod_{l=1}^N \binom{x_l}{n_l - m_l}$$

The more useful proposition for computing can be deduced from this theorem easily

**Proposition 6.4.1.** The difference operator is the product of the partial difference operators. i.e

$$\Delta^{\underline{m}} \prod_{l=1}^N \binom{x_l}{n_l} = \prod_{l=1}^N \Delta_l^{m_l} \binom{x_l}{n_l}$$

The following somewhat says that it is meaningful to consider $\mathbb{B}$−polyfract as a function, i.e. the function with $\text{fdeg}(f) < \infty$ is $\mathbb{B}$−ployfract and the definition is introduced formally in last section.

**Theorem 6.5** ((3), pg.9). [The Fundamental Representation Theorem] Let $f: \mathbb{Z}^N \to \mathbb{B}$ be a map with $\text{fdeg}(f) < \infty$. Then

$$f(\underline{x}) = \sum_{\underline{n} \in \mathcal{N}} \binom{x_1}{n_1}\cdots\binom{x_N}{n_N} b_{\underline{n}}$$

for some finite $\mathcal{N} \subset \mathbb{N}^N$ and some finite $\mathcal{A} \subset \mathbb{B}$ satisfying $\text{fdeg}(f) = \max\{|\underline{n}| \mid \underline{n} \in \mathcal{N}\}$.

**Theorem 6.6.** Consider a $\mathbb{B}$−polyfract $(\mathcal{P}, \{b\}, \{\underline{n}\})$. Then

$$\text{count}(\mathcal{P}) = \text{fdeg}(f(\mathcal{P}))$$

*Proof.* By Theorem 6.5: we can write

$$f(\mathcal{P})(\underline{n}) = \binom{x_1}{n_1} \cdots \binom{x_N}{n_N} b$$

By definition: $\text{count}(\mathcal{P}) = |\underline{n}|$. Applying Lemma 6.2 and Lemma 6.3 to $f(\mathcal{P})$ with $\mathbb{A} = \mathbb{Z}^N$ and $\mathbb{B} = \mathbb{Z}$, it is enough to show $\Delta_{e_1}^{|\underline{n}|} f(\mathcal{P})(\underline{0}) \neq 0$ and there exists some $1 \leq i \leq N$ such that $\Delta_{e_i}^{|\underline{n}|+1} f(\mathcal{P}) \equiv 0$. We are going to finish the proof by computing. To use Theorem 6.4, denoting $\underline{m} = (|\underline{n}|, 0, \cdots, 0)$: we have

$$\Delta^{\underline{m}} = \Delta_1^{|\underline{n}|} = \Delta_{e_1}^{|\underline{n}|}$$

and hence

$$\Delta_{e_1}^{|\underline{n}|} f(\mathcal{P})(\underline{x}) = \Delta^{\underline{m}} f(\mathcal{P})(\underline{x})$$
$$= \prod_{l=1}^{N} \binom{x_l}{n_l - m_l} b$$
$$= \binom{x_1}{n_1 - |\underline{m}|} \prod_{l=2}^{N} \binom{x_l}{n_l} b$$

Now, for $1 \leq l \leq N$, denote $\text{bi}_l(x, y) = \binom{x}{y}$, then by Lemma 6.1

$$\Delta_{e_1}^{|\underline{n}|} f(\underline{x}) = \text{bi}_1(x_1, n_1 - |\underline{n}|) \prod_{l=2}^{N} \text{bi}_2(x_l, n_l)$$

We evaluate this at $\underline{0}$:

$$\Delta_{e_1}^{\underline{n}} f(\mathcal{P})(\underline{0}) = \text{bi}_1(0, n_1 - |\underline{n}|) \prod_{l=2}^{N} \text{bi}_2(0, n_l)$$

Now since $\underline{n} \in \mathbb{N}_{\geq 1}$ implying for each $2 \leq l \leq N$, $n_l \geq 1$, applying the first item of Lemma 6.1 to $2 \leq l \leq N$:

$$\prod_{l=2}^{N} \text{bi}_2(0, n_l) = 0$$

Also we have $n_1 - |\underline{n}| < 0$, applying the fourth item of Lemma 6.1 to $\text{bi}_1$, $\text{bi}_1(0, n_1 - |\underline{n}|) = 0$. And hence their product is zero. To find such $i$ and verify it, we choose $\underline{p} \in \mathbb{N}^N$ such that $|\underline{p}| = |\underline{n}| + 1$. Then $p_1 + \cdots + p_N = n_1 + \cdots + n_N + 1$ implies there must be $i$ such that $p_i \geq n_i + 1$. denote $\underline{\omega} = \underbrace{(0, \cdots, 0, p_i, 0, \cdots, 0)}_{p_i \text{ is in } i\text{th position}}$. We have by our choice of $\underline{p}$

$$\Delta_{e_i}^{|\underline{n}|+1} = \Delta_{e_i}^{|\underline{p}|} = \Delta_{e_i}^{|\underline{p}| - p_i} \Delta_{e_i}^{p_i}$$

Now to verify such $i$: it is enough to prove $\Delta_{e_i}^{p_i} f(\mathcal{P}) \equiv 0$, since the derivative of the zero map is still zero. Since

$$\Delta^{\underline{\omega}} = \Delta_i^{p_i} = \Delta_{e_i}^{p_i}$$

Then denoting $\text{bi}_i(x_i, y_i) = \binom{x_i}{y_i}$, by Theorem 6.4: for each $\underline{x} \in \mathbb{Z}^N$ by the commutative

$$\Delta_{e_i}^{p_i} f(\mathcal{P})(\underline{x}) = \Delta^{\underline{\omega}} f(\mathcal{P})(\underline{x})$$
$$= \prod_{l=1}^{N} \binom{x_l}{n_l - \omega_l} b$$
$$= \binom{x_i}{n_i - \omega_i} \prod_{l=1, l \neq i}^{N} \binom{x_l}{n_l - \omega_l} b$$
$$= \binom{x_i}{n_i - p_i} \prod_{l=1, l \neq i}^{N} \binom{x_l}{n_l - \omega_l} b$$
$$= \text{bi}_i(x_i, n_i - p_i) \prod_{l=1, l \neq i}^{N} \binom{x_l}{n_l - \omega_l} b$$

since $p_i \geq n_i + 1$ implying $n_i - p_i < 0$, applying the fourth item of Lemma 6.1 to $\text{bi}_i$, we have $\text{bi}_i(x_i, n_i - p_i) = 0$. This implies $\Delta_{e_i}^{p_i} f(\mathcal{P}) \equiv 0$ ∎

**Theorem 6.7.** Consider a multivariate polynomial function $f : \mathbb{Z}^N \to \mathbb{Z}$ defined by
$$f(\underline{x}) = \sum_{\underline{n} \in \mathcal{N}} \binom{x_1}{n_1} \cdots \binom{x_N}{n_N} b_{\underline{n}}$$
where $\mathcal{N} \subset \mathbb{N}^N$ with $\operatorname{card}(\mathcal{N}) < \infty$. Let $\underline{m} \in \mathcal{N}$ be one satisfying $|\underline{m}| = \max\{|\underline{n}| \mid \underline{n} \in \mathcal{N}\}$. Denote $g : \mathbb{Z}^N \to \mathbb{B}$ by
$$g(\underline{x}) = \binom{x_1}{m_1} \cdots \binom{x_N}{m_N} b_{\underline{n}}$$
Then $\operatorname{fdeg}(f) = \operatorname{fdeg}(g)$

*Proof.* After combining those exactly same term monomials, we can write $\mathcal{N} = \{\underline{n}_1, \underline{n}_2, \cdots, \underline{n}_k\}$ where $k := \operatorname{card}(\mathcal{N})$ and without loss of generality we can assume that $|\underline{n}_1| < |\underline{n}_2| < |\underline{n}_3| < \cdots < |\underline{n}_k|$ and hence making sure such $\underline{m}$ is the only one satisfying the condition stated in the statement. Construct $\mathbb{Z}$–polyfract $(\mathcal{P}, \{b_{\underline{m}}\}, \{\underline{m}\})$ and applying Theorem 6.6 to it, we have $\operatorname{fdeg}(g) = \operatorname{count}(\mathcal{P}(\{b_{\underline{m}}\}, \{\underline{m}\})) = |\underline{m}|$. And hence $\operatorname{count}(\mathcal{P}) = |\underline{n}_k| = |\underline{m}|$. For each $\underline{n}_i$, denote $g_{\underline{n}_i}(\underline{x}) = \binom{x_1}{n_i^{(1)}} \cdots \binom{x_N}{n_i^{(N)}} b_{\underline{n}_i}$. Then we have by the condition on $\underline{m}$

$$f(\underline{x}) = \sum_{\underline{n} \in \mathcal{N}} \binom{x_1}{n_1} \cdots \binom{x_N}{n_N} b_{\underline{n}}$$
$$= \sum_{l=1}^{k} \binom{x_1}{n_l^{(1)}} \cdots \binom{x_N}{n_l^{(N)}} b_{\underline{n}_l}$$
$$= \sum_{l=1}^{k} g_{\underline{n}_i}(\underline{x})$$

Now compute that by the linearity of the polynomial functions, i.e. ♣ in Lemma 6.3,

$$\Delta^{\underline{m}} f = \Delta^{\underline{m}} \sum_{l=1}^{k} g_{\underline{n}_l}$$
$$= \sum_{l=1}^{k} \Delta^{\underline{m}} g_{\underline{n}_l}$$
$$= \sum_{l=1}^{k-1} \Delta^{\underline{m}} g_{\underline{n}_l} + \Delta^{\underline{m}} g_{\underline{m}} \tag{0.15}$$

Now we are going to say that for all $1 \leq l \leq k-1$, $\Delta^{\underline{m}} g_{\underline{n}_l} = 0$ and hence $\Delta^{\underline{m}} f = \Delta^{\underline{m}} g_{\underline{m}}$. We only consider the case $l = k-1$ and other cases can be explained successively. Compute by Theorem 6.4:

$$\Delta^{\underline{m}} g_{\underline{n}_{k-1}}(\underline{x}) = \Delta^{\underline{m}} \binom{x_1}{n_{k-1}^{(1)}} \cdots \binom{x_N}{n_{k-1}^{(N)}} b_{\underline{n}_{k-1}}$$
$$= \binom{x_1}{n_{k-1}^{(1)} - m_1} \cdots \binom{x_N}{n_{k-1}^{(N)} - m_N} b_{\underline{n}_{k-1}}$$

Since $\underline{n}_{k-1} < \underline{n}_k = \underline{m}$, there must be some coordinate $n_{k-1}^{(p)}$ of $\underline{n}_{k-1}$ such that $n_{k-1}^{(p)} < m_p = n_k^{(p)}$ and hence $\binom{x_p}{n_{k-1}^{(p)} - m_p} \equiv 0$ and by the production to conclude $\Delta^{\underline{m}} g_{\underline{n}_{k-1}} = 0$. Now by our reordering on $\mathcal{N}$, we have $\Delta^{\underline{m}} f = \Delta^{\underline{m}} g_{\underline{m}} = \Delta^{\underline{m}} g \neq 0$ and indeed among the product of operators of $\Delta^{\underline{m}} = \Delta_{e_N}^{m_N} \cdots \Delta_{e_1}^{m_1}$, there are $|\underline{m}|$ elements chosen from $\mathbb{Z}^N$. To finish the proof, it remains to show there exists $\underline{m}+1$ elements of $\mathbb{Z}^N$ to form discrete difference operators so that after taking operators of $f$ under this difference operators, the result polynomial map is identical to zero. By the identity 0.15 and definition of maps, it is enough to find such operators for $g_{\underline{m}}$, i.e. $g$. Now, consider $\underline{m}' := (m_1, m_2, \cdots, m_N + 1)$ and $\Delta^{\underline{m}'} g$ and by the similarly discussion, it enough to show $\Delta^{\underline{m}'} g = 0$. By Theorem 6.4:

$$\Delta^{\underline{m}'} g(\underline{x}) = \Delta^{\underline{m}'} g_{\underline{m}}(\underline{x})$$
$$= \binom{x_1}{n_k^{(1)} - m_1} \cdots \binom{x_N}{n_k^{(N)} - m_N - 1} b_{\underline{m}}$$
$$= 0 \cdot b_{\underline{m}}$$
$$= b_{\underline{m}}$$

where the second step from the last is by the uniqueness, i.e $\underline{m} \in \mathcal{N}$ such that $|\underline{m}| = |\underline{n}_k|$. Then $n_k^N = m_N$ and hence $\binom{x_N}{n_k^{(N)} - m_N - 1} \equiv \binom{x_N}{-1} \equiv 0$ and the production and the way we forms $\mathbb{Z}$–polyfract polynomial function is by the group action $\operatorname{End}(\mathbb{Z}^{\mathbb{Z}^N}) \times \mathbb{Z}^{\mathbb{Z}^N} \to \mathbb{Z}^{\mathbb{Z}^N}$.

∎

Now we present the main theorem and are ready to prove it.

**Theorem 6.8.** The counting number of a $\mathbb{Z}$−polyfract $\mathcal{P}$ is determined by the number of taking discrete difference operators of its associated function. Let a triple $(\mathcal{P}, \mathcal{A}, \mathcal{N})$ be a $\mathbb{Z}$−polyfract, i.e some finite $\mathcal{A} \subset \mathbb{Z}$ and some finite $\mathcal{N} \subset \mathbb{N}^N$ Then

$$\text{count}(\mathcal{P}) = \text{fdeg}(f(\mathcal{P}))$$

*Proof.* Denote $k := \text{card}(\mathcal{N})$. Write

$$\mathcal{N} = \{\underline{n}_1, \underline{n}_2, \cdots, \underline{n}_k\}$$

Now reordering it (or without of loss of generality, assume it) as

$$|\underline{n}_1| \leq |\underline{n}_2| \leq |\underline{n}_3| \leq \cdots \leq |\underline{n}_k|$$

Also if the increasing number between consequences is strictly bigger than one, we can always insert a monomial term. In detail, there exists $l \in \{1, 2, \cdots, k\}$ such that $|\underline{n}_{l+1}| - |\underline{n}_l| \geq 2$. Denote $\underline{n}_{l+1} = (n_{l+1}^{(1)}, n_{l+1}^{(2)}, \cdots, n_{l+1}^{(N)})$ and $\underline{n}_l = (n_l^{(1)}, n_l^{(2)}, \cdots, n_l^{(N)})$. Denote $\mathcal{N}_l = \{\underline{n}_l\}$, $\mathcal{N}_{l+1} = \{\underline{n}_{l+1}\}$ and $\mathcal{A}_l = \{b_{\underline{n}} \mid \underline{n} \in \mathcal{N}_l\}$, $\mathcal{A}_{l+1} = \{b_{\underline{n}} \mid \underline{n} \in \mathcal{N}_{l+1}\}$. Now write

$$\mathcal{P}(\mathcal{A}_l, \mathcal{N}_l) = \binom{x_1}{n_l^{(1)}} \cdots \binom{x_N}{n_l^{(N)}} b_{\underline{n}_l} \quad \mathcal{P}(\mathcal{A}_{l+1}, \mathcal{N}_{l+1}) = \binom{x_1}{n_{l+1}^{(1)}} \cdots \binom{x_N}{n_{l+1}^{(N+1)}} b_{\underline{n}_{l+1}}$$

Now we only consider the case: $|\underline{n}_{l+1}| - |\underline{n}_l| = 2$ and $n_{l+1}^{(1)} - n_l^{(1)} = 2$ and other cases can be discussed similarly. We can union $\mathcal{N}$ with a subset and then

$$\binom{x_1}{n_l^{(1)} + 1} \cdots \binom{x_N}{n_l^{(N)}} b_{\underline{n}_l}$$

For each $1 \leq l \leq k$ assume $|\underline{n}_{l+1}| - |\underline{n}_l| < 2$ and construct $k$ number of $\mathbb{Z}$−polyfract triples as

$$\mathcal{N}_1 = \{\underline{n}_1\}, \mathcal{N}_2 = \{\underline{n}_2\}, \cdots, \mathcal{N}_k = \{\underline{n}_k\}$$
$$\mathcal{A}_1 = \{b_{\underline{n}} \mid \underline{n} \in \mathcal{N}_1\}, \mathcal{A}_2 = \{b_{\underline{n}} \mid \underline{n} \in \mathcal{N}_2\}, \cdots, \mathcal{A}_k = \{b_{\underline{n}} \mid \underline{n} \in \mathcal{N}_k\}$$

Then

$$\text{count}(\mathcal{P}(\mathcal{A}_1, \mathcal{N}_1)) \leq \text{count}(\mathcal{P}(\mathcal{A}_2, \mathcal{N}_2)) \cdots \leq \text{count}(\mathcal{P}(\mathcal{A}_k, \mathcal{N}_k)) = \text{count}(\mathcal{P}) \tag{0.16}$$

which says we only need to consider the sub-collection $\mathcal{P}(\mathcal{A}_k, \mathcal{N}_k)$. Denote

$$\underline{n}_k = (n_k^{(1)}, n_k^{(2)}, \cdots, n_k^{(N)})$$

Then the associated polynomial function $f(\mathcal{P}(\mathcal{A}_k, \mathcal{N}_k))(\underline{x})$ is of the form

$$f(\mathcal{P}(\mathcal{A}_k, \mathcal{N}_k))(x_1, x_2, \cdots, x_N) = \binom{x_1}{n_k^{(1)}} \cdots \binom{x_N}{n_k^{(N)}} b_{\underline{n}_k}$$

Now by Theorem 6.7,

$$\text{fdeg}(f) = \text{fdeg}(f(\mathcal{P}(\mathcal{A}_k, \mathcal{N}_k))) \tag{0.17}$$

and by Theorem 6.6

$$\text{count}(\mathcal{P}(\mathcal{A}_k, \mathcal{N}_k)) = \text{fdeg}(f(\mathcal{P}(\mathcal{A}_k, \mathcal{N}_k)) \tag{0.18}$$

Then combining the equality 0.16, the equality 0.18 and the equality 0.17, the result is achieved. ∎

Since all computing techniques used in the left hand side of the equality are standard basis ones, i.e $\{\Delta_{e_1}, \cdots, \Delta_{e_N}\}$ and all used in the right hand side of it are arbitrary ones, i.e $\{\Delta_a \mid a \in \mathbb{Z}^N\}$, we actually finish the progress of linearization. But to make our achievement more original one in mathematical statement, we introduce the following notations and recall some more fundamental notations for presenting our statement. We introduce the following notation $\text{fdeg}^{\mathbb{A}}(f) = \max\{d \in \mathbb{N}_{\geq 1} \mid \Delta_{a_1} \Delta_{a_2} \cdots \Delta_{a_d} f \neq 0 \text{ and } (a_1, \cdots, a_d) \in \mathbb{A}^d\}$. Since the value of fdeg is exactly the same with the value of the usual degree, with a little contrapositive deduction, $\text{fdeg}^{\mathbb{A}}(f) = \min\{d \in \mathbb{N}_{\geq 1} \mid \Delta_{a_1} \cdots \Delta_{a_{d+1}} f = 0 \text{ and } (a_1, \cdots, a_{d+1}) \in \mathbb{A}^{d+1}\}$ and recall

$$\Delta_a f(x) := f(x + a) - f(x)$$

we define the degree of the constant function to be 1 and we are not interested in the trivial case $f = 0$. We also introduce $\text{fdeg}^{\mathbb{Z}^N}(f)$. We also introduce the following notation $\text{fdeg}^{\mathbb{Z}^N}(f) = \min\{|\underline{n}| - 1 \mid \Delta^{\underline{n}} f = 0 \text{ and } \underline{n} \in \mathbb{N}_{\geq 1}^N\}$, where $\Delta^{\underline{n}} := \Delta_1^{n_1} \cdots \Delta_N^{n_N}$ and $\Delta_k^k$ is $k$−fold product of $\Delta_k$ for each $1 \leq k \leq N$. Similarly $\text{fdeg}^{\mathbb{Z}^N}(f) = \max\{|\underline{n}| \in \mathbb{N}_{\geq 1} \mid \Delta^{\underline{n}} f \neq 0 \text{ and } \underline{n} \in \mathbb{N}_{\geq 1}^N\}$, where $|\underline{n}| := (n_1, \cdots, n_N)$. By analogue in $\mathbb{R}^N$, we call $\Delta_{a_1} \cdots \Delta_{a_d}$ as multivariate discrete difference operators. Now with the help of Theorem 6.8, we can have the following theorem

**Theorem 6.9.** Let $f : \mathbb{Z}^N \to \mathbb{Z}$ be a map of finite functional degree. Then

$$\mathrm{fdeg}^{\mathbb{A}}(f) = \mathrm{fdeg}^{\mathbb{Z}^N}(f)$$

We must emphasize that there is no difference of the distinction between polynomial functions and the associated function of polynomials.

# 7 Appendix

*Proof of Theorem 3.4.* The first expansion:
We use the following process to prove this formula: For each bracket, we use the same index dummy variable for $n$ pairs of brackets in the summation way:

$$(x_{1,1} + x_{1,2} + x_{1,3} + \cdots + x_{1,m_1})(x_{2,1} + x_{2,2} + x_{2,3} + \cdots + x_{2,m_2}) \cdots (x_{n,1} + x_{n,2} + x_{n,3} + \cdots + x_{n,m_n})$$

$$= \sum_{l=1}^{m_1} x_{1,l} \sum_{l=1}^{m_2} x_{2,l} \cdots \sum_{l=1}^{m_n} x_{n,l}$$

For the $n$ times product, we enumerate the index variables to $p_1, p_2, p_3 \cdots, p_n$:

$$\sum_{l=1}^{m_1} x_{1,l} \sum_{l=1}^{m_2} x_{2,l} \cdots \sum_{l=1}^{m_n} x_{n,l} = \sum_{p_1=1}^{m_1} x_{1,p_1} \sum_{p_2=1}^{m_2} x_{2,p_2} \cdots \sum_{p_n=1}^{m_n} x_{n,p_n}$$

Applying distributive law $n$ times from the outside to inside

$$\sum_{p_1=1}^{m_1} x_{1,p_1} \sum_{p_2=1}^{m_2} x_{2,p_2} \sum_{p_3=1}^{m_3} x_{3,p_3} \cdots \sum_{p_n=1}^{m_n} x_{n,p_n} = \sum_{p_1=1}^{m_1} \sum_{p_2=1}^{m_2} x_{1,p_1} x_{2,p_2} (\sum_{p_3}^{m_3} x_{3,p_3} \cdots \sum_{p_n=1}^{m_n} x_{n,p_n})$$

$$\vdots$$

$$= \sum_{p_1=1}^{m_1} \sum_{p_2=1}^{m_2} \cdots \sum_{p_n=1}^{m_n} x_{1,p_1} x_{2,p_2} \cdots x_{n,p_n}$$

The second expansion: By the distributive law, the $n$th degree monomial term in the final summation is of the form

$$x_{1,p_1} x_{2,p_2} x_{3,p_3} \cdots x_{n,p_n}$$

To find index under summation, observing that $p_1 \in \{1, 2, 3, \cdots, m_1\}$, $p_2 \in \{1, 2, 3, \cdots, m_2\}$, $\cdots$, $p_n \in \{1, 2, 3, \cdots, m_n\}$, according to the distributive law, taking product of sets: $(p_1, p_2, \cdots, p_n) \in [m_1] \times [m_2] \times \cdots \times [m_n]$, denoting $[m_i] := \{1, 2, 3, \cdots, m_i\}$ for each $1 \leq i \leq n$ and putting it under the summation:

$$\sum_{(p_1,p_2,\cdots,p_n)\in[m_1]\times[m_2]\times\cdots\times[m_n]} x_{1,p_1} x_{2,p_2} \cdots x_{n,p_n}$$

Using $\prod$ for the index and the terms:

$$\sum_{(p_1,p_2,\cdots,p_n)\in\prod_{i=1}^{n}[m_i]} \left\{ \prod_{i=1}^{n} x_{i,p_i} \right\}$$

∎

**Theorem 7.1** ((17), pg.19). For real $n$:

$$\binom{-n}{k} = (-1)^k \binom{n+k-1}{k}$$

*Proof.* Just multiplying each factor of denominators by $-1$ one-by-one and reverse them. ∎

**Theorem 7.2.** (3, pg.7) Let $\mathbb{A}$ and $\mathbb{B}$ be commutative groups. Let $a \in \mathbb{A}$. Let $n \in \mathbb{N}_{\geq 2}$. Let $\Delta_a^n$ be $n$−fold product $\Delta_a \cdots \Delta_a \in \mathrm{End}(\mathbb{B}^{\mathbb{A}})$. Let $f \in \mathbb{B}^{\mathbb{A}}$ and $x \in \mathbb{A}$. Then

$$\Delta_a^n f(x) = \sum_{i=0}^{n} (-1)^i \binom{n}{i} f(x + (n-i)a) = \sum_{i=0}^{n} (-1)^{n-i} \binom{n}{i} f(x + ia)$$

**Theorem 7.3.** Let $h: \mathbb{Z}^N \to \mathbb{Z}$ be a multivariate polynomial such that its corresponding $\mathbb{Z}$−polyfract $\mathcal{P}(\mathcal{A}, \mathcal{N})$ satisfying $\mathcal{A} = \{1\} \subset \mathbb{Z}$ and $\mathcal{N} \subset \mathbb{N}_{\geq 1}^N$ with card$(\mathcal{N}) = 1$. i.e. $h = f(\mathcal{P})$. Denoting $\mathcal{N} = \{\underline{m}\} \subset \mathbb{N}_{\geq 3}^N$. Let $\underline{n} \in \mathbb{N}_{\geq 2}^N$ such $n_i < m_i$ for each $1 \leq i \leq N$. Then for each $\underline{x} = (x_1, \cdots, x_N) \in \mathbb{Z}^N$:

$$\sum_{0 \leq p_1 \leq n_1; 0 \leq p_2 \leq n_2; \cdots; 0 \leq p_N \leq n_N} (-1)^{p_1 + \cdots + p_N} h(\sum_{i=1}^N (n_i - p_i + x_i) e_i) = \prod_{l=1}^N \binom{x_l}{m_l - n_l}$$

*Proof.* For proving it, since $\mathbb{Z}^N = \langle e_1, \cdots, e_N \rangle_{\mathbb{Z}}$, by using Theorem 7.2 to $h$ repeatedly to each $e_i$, computing: applying it to $\Delta_{e_1}^{n_1}$: for $\underline{x} \in \mathbb{Z}^N$

$$\Delta_{e_1}^{n_1} h(\underline{x}) = \sum_{p_1=0}^{n_1} (-1)^{p_1} \binom{n_1}{p_1} h(\underline{x} + (n_1 - p_1) e_1)$$

Denoting $g(\underline{x}) := \Delta_{e_1}^{n_1} h(\underline{x})$, applying it again to $g(\underline{x})$:

$$\begin{aligned}
\Delta_{e_2}^{n_2} \Delta_{e_1}^{n_1} h(\underline{x}) &= \Delta_{e_2}^{n_2} g(\underline{x}) \\
&= \sum_{p_2=0}^{n_2} (-1)^{p_2} \binom{n_2}{p_2} g(\underline{x} + (n_2 - p_2) e_2) \\
&= \sum_{p_2=0}^{n_2} (-1)^{p_2} \binom{n_2}{p_2} \Delta_{e_1}^{n_1} h(\underline{x} + (n_2 - p_2) e_2) \\
&= \sum_{p_2=0}^{n_2} (-1)^{p_2} \binom{n_2}{p_2} \sum_{p_1=0}^{n_1} (-1)^{p_1} \binom{n_1}{p_1} h(\underline{x} + (n_2 - p_2) e_2 + (n_1 - p_1) e_1) \\
&= \sum_{p_2=0}^{n_2} \sum_{p_1=0}^{n_1} (-1)^{p_2} (-1)^{p_1} h(\underline{x} + (n_2 - p_2) e_2 + (n_1 - p_1) e_1)
\end{aligned}$$

Inductively it is reasonable to deduce

$$\begin{aligned}
\text{R.H.S} &:= \Delta_{e_N}^{n_N} \cdots \Delta_{e_2}^{n_2} \Delta_{e_1}^{n_1} h(\underline{x}) \\
&= \sum_{p_N=0}^{n_N} \cdots \sum_{p_1=0}^{n_1} (-1)^{p_N} (-1)^{p_{N-1}} \cdots (-1)^{p_1} h(\underline{x} + (n_N - p_N) e_N + \cdots + (n_1 - p_1) e_1) \\
&= \sum_{0 \leq p_1 \leq n_1; 0 \leq p_2 \leq n_2; \cdots; 0 \leq p_N \leq n_N} (-1)^{p_1 + \cdots + p_N} h(\underline{x} + \sum_{i=1}^N (n_i - p_i) e_i) \\
&= \sum_{0 \leq p_1 \leq n_1; 0 \leq p_2 \leq n_2; \cdots; 0 \leq p_N \leq n_N} (-1)^{p_1 + \cdots + p_N} h(\sum_{i=1}^N (n_i - p_i + x_i) e_i)
\end{aligned}$$

Denoting $\underline{n} = (n_1, \cdots, n_N)$, Use Theorem 6.4 computing:

$$\begin{aligned}
\text{L.H.S} &:= \Delta_{e_N}^{n_N} \cdots \Delta_{e_1}^{n_1} h(\underline{x}) \\
&= \Delta^{\underline{n}} h(\underline{x}) \\
&= \prod_{l=1}^N \binom{x_l}{m_l - n_l}
\end{aligned}$$

Finally since R.H.S = L.H.S, the equality is achieved. ∎